# Elementary Primes Counting Methods

**N. A. Carella, August 2012.**

***Abstract*** This work proposes elementary proofs of several related primes counting problems based on an elementary weighted sieve. The subsets of primes considered here are the followings: the subset of twin primes $P_T = \{ p$ and $p + 2$ are primes $\}$, the subset of Germain primes $P_G = \{ p$ and $2p + 1$ are primes $\}$, and the subset of quadratic primes $P_f = \{ p = n^2 + 1$ primes $\}$. These subsets of primes are widely believed to be infinite subsets of prime numbers.

**Table of Content**

# 1. Introduction

The problem of determining the cardinality of the subset of primes $P_f = \{ p = f(n) : p \text{ is prime and } n \geq 1 \}$, defined by a function $f : \mathbb{N} \to \mathbb{N}$, as either finite or infinite is vastly simpler than the problem of determining the primes counting function of the subset of primes $\pi_f(x) = \#\{ p = f(n) \leq x : p \text{ is prime and } n \geq 1 \}$.

The cardinality of the simplest case of the set of all primes $\mathbb{P} = \{ 2, 3, 5, \ldots \}$ was settled by Euclid over two millennia ago, see [EU]. After that event, many other proofs have been discovered by other authors, see [RN, Chapter 1], [NW, p. 4], and the literature. In contrast, the determination of the primes counting function of the set of all primes $\pi(x) = \#\{ p \leq x : p \text{ is prime} \}$ is still incomplete. Indeed, the current asymptotic formula (99) of the primes counting function, also known as the Prime Number Theorem, is essentially the same as determined by delaVallee Poussin, and Hadamard about a century ago.

The main results are the possible proofs, based on elementary methods, of the followings primes counting problems:

1. The subset of twin primes $P_T = \{ p \text{ and } p + 2 \text{ are primes} \}$, Theorem 2.1,
2. The subset of Germain primes $P_G = \{ p \text{ and } 2p + 1 \text{ are primes} \}$, Theorem 3.1,
3. The subset of quadratic primes $P_f = \{ p = n^2 + 1 \text{ primes} \}$, Theorem 4.1.

These subsets of primes are widely believed to be infinite subsets of prime numbers, and confirmed here. Other results included are a completely elementary proof of the prime number theorem on arithmetic progressions in Theorem 5.1, and a completely elementary proof of the prime number theorem in Theorem 7.2.

The techniques used on these problems are applicable to some other classes of primes, but these techniques probably do work on subsets of primes defined by exponential functions. For example, the subset of Fermat numbers $\{ F_n = 2^{2^n} + 1 : n \geq 1 \}$, and the subset of Mersenne numbers $\{ M_p = 2^p - 1 : p \geq 2 \text{ prime} \}$.

# 2. Twin Primes

The twin primes conjecture claims that the Diophantine equation $y = x + 2$ has infinitely many prime pairs solutions. More generally, the dePolignac conjecture claims that for any fixed $k \geq 1$, the Diophantine equation $y = x + 2k$ has infinitely many prime pairs solutions, confer [RN, p. 265], [GS], [KV], [FI, p. 315], [AR], [KR], [NW, p. 337], [WK], [WS], [PP], and related topics in [BC], [GP], [GU, p. 31], [SD], et alii. This work proposes a proof of this primes counting problem. In particular, the cardinality of the subset of twin primes $\{ p \text{ and } p + 2 \text{ are primes} \}$ is infinite.

### 2.1 The Main Result

***Theorem* 2.1.** For any fixed $k \geq 1$, the Diophantine equation $y = x + 2k$ has infinitely many prime pairs $x = p$, $y = q$ solutions.

Proof: Without loss in generality, let $k = 1$, and consider the weighted finite sum over the integers:

$$\begin{aligned}\sum_{n \le x} n\Lambda(n)\Lambda(n+2) &= -\sum_{n \le x} n\Lambda(n) \sum_{d \mid n+2} \mu(d)\log d \\ &= -\sum_{d \le x+2} \mu(d)\log d \left( \sum_{n \le x,\ n \equiv -2 \bmod d} n\Lambda(n) \right),\end{aligned} \tag{1}$$

where $\Lambda(n) = -\sum_{d \mid n} \mu(d)\log d$. This follows from Lemma 6.1 in Section 6, and inverting the order of summation.

Other examples of inverting the order of summation are given in [MV, p. 35], [RE, p. 27], [RM, p. 216], [SP, p. 83], [TM, p. 36], et alii. The proof of Theorem 7.3 is quite similar to this, and it is instructive to review it.

Applying Lemmas 6.2, and 6.3, (with $a = -2$, and $q = d$), and simplifying yield

$$\begin{aligned}\sum_{n \le x} n\Lambda(n)\Lambda(n+2) &= -\sum_{d \le x+2,\ d \text{ odd}} \mu(d)\log d \left( \frac{1}{2\varphi(d)} x^2 + \Omega_\pm\left(x^{3/2}\log_3 x\right)\right) \\ &\ge -\sum_{d \le x+2,\ d \text{ odd}} \mu(d)\log d \left( \frac{1}{2d} x^2 \right)\end{aligned} \tag{2}$$

for infinitely many real numbers $x \ge 1$. The restriction to odd $d \ge 1$ stems from the relation $\Lambda(n) = 0$ for even $d \ge 2$. This restriction is equivalent to $\gcd(a, q) = 1$ with $a = -2$, and $q = d$ in Lemmas 6.2 and 6.3. Here, there is no restriction on the range of $d \le x$.

Applying Lemma 6.4, (for example, at $s = 1$, and $q_0 = 4$), the previous equation becomes

$$\begin{aligned}\sum_{n \le x} n\Lambda(n)\Lambda(n+2) &\ge -c_1 x^2 \sum_{d \le x+2,\ d \text{ odd}} \frac{\mu(d)\log d}{d} \\ &= c_1 x^2 \left( c_2 + O\left( e^{-c(\log x)^{1/2}} \log x \right)\right) \\ &= c_3 x^2 + O\left( x^2 e^{-c(\log x)^{1/2}} \log x \right),\end{aligned} \tag{3}$$

where $c_1, c_2, c_3 > 0$ and $c > 0$ are nonnegative constants, for infinitely many real numbers $x \ge 1$.

Now, assume that there are finitely many twin primes $p$, and $p + 2 < x_0$, where $x_0 > 0$ is a large constant. For example, $\Lambda(n)\Lambda(n + 2) = 0$ for all $n > x_0$. Then

$$c_3 x^2 + O\left( x^2 e^{-c(\log x)^{1/2}} \log x \right) \le \sum_{n \le x} n\Lambda(n)\Lambda(n+2) \le x_0^3 \tag{4}$$

for all real numbers $x \ge 1$. But this is a contradiction for infinitely many real numbers $x > x_0^2$. Ergo, there are infinitely many intervals $[1, x]$ that contain

$$\pi_2(x) \ge c_3 \frac{x}{\log^2 x} + O\left( x e^{-c(\log x)^{1/2}} \right) \tag{5}$$

twin primes, see (9), as $x \geq 1$ tends to infinity. Nil volentibus arduum. ■

Another way to derive this result is to use Lemma 6.2-ii, or 6.2-iii. Proceed to break the finite sum (1) into two finite sums:

$$\sum_{n \leq x} n\Lambda(n)\Lambda(n+2) = -\sum_{d \leq M} \mu(d)\log d \sum_{\substack{n \leq x/d, \\ n \equiv -2 \bmod d}} n\Lambda(n) - \sum_{d > M} \mu(d)\log d \sum_{\substack{n \leq x/d, \\ n \equiv -2 \bmod d}} n\Lambda(n) \tag{6}$$

where $M = O(\log x)^B$, $B > 0$ constant. And determine appropriate estimates for these finite sums. But, the one given above is much simpler, see [MV, p. 386] for the analysis of similar counting problems. It should also be noted that other (twin primes sieve) weighted finite sums over the integers such as

$$\sum_{n \leq x} \varphi(n)\Lambda(n)\Lambda(n+2), \qquad \sum_{n \leq x} \sigma(n)\Lambda(n)\Lambda(n+2), \qquad \sum_{n \leq x} \tau(n)\Lambda(n)\Lambda(n+2), \tag{7}$$

where φ, σ, and τ are the totient function, and the divisors functions respectively, yield precisely the same result. The analysis of these later finite sums are quite similar to the Titchmarsh divisor problem over arithmetic progressions, see [CM, p. 172] for an introduction to this problem.

## 2.2 Counting Functions

Previous works on the twin primes conjecture have studied certain Dirichlet series, see [AR], [KR], et cetera. The primes counting method introduced here suggests that a more suitable Dirichlet series is the complex valued function

$$D_{2k}(s) = \sum_{n \geq 1} \frac{n\Lambda(n)\Lambda(n+2k)}{n^s} = \int_1^\infty \frac{dR(x)}{x^s}, \tag{8}$$

where $R(x) = \sum_{n \leq x} n\Lambda(n)\Lambda(n+2)$, $k \geq 1$ is a fixed integer, and $\Re e(s) \geq 2$. This series leads to another approach to study of the twin primes conjecture, and more generally, the dePolignac conjecture via the Wiener-Ikehara theorem introduced in [AR].

An application of partial summation to the finite sum yields the counting function

$$\begin{aligned} \pi_2(x) &= \sum_{n \leq x} \frac{n\Lambda(n)\Lambda(n+2)}{n\log(n)\log(n+2)} - \sum_{n=p^k \text{ or } n+2=q^m \leq x,\ k,m \geq 2} \frac{n\Lambda(n)\Lambda(n+2)}{n\log(n)\log(n+2)} \\ &= \int_2^x \frac{dR(t)}{t\log(t)\log(t+2)} + O(x^{1/2}\log^2 x), \end{aligned} \tag{9}$$

where $R(x) = \sum_{n \leq x} n\Lambda(n)\Lambda(n+2)$. Evaluating the integral realizes a lower estimate of the twin primes counting function

$$\pi_2(x) = \#\{\, p \leq x : p \text{ and } p+2 \text{ are primes} \,\} \geq c_1 \frac{x}{\log^2 x} + O(\frac{x}{\log^3 x}). \tag{10}$$

This lower estimate has the correct order of magnitude as determined by sieve method, that is,

$$\pi_2(x) = \#\{\ p \le x : p \ \text{ and } \ p+2 \ \text{ are primes}\ \} \le c_3 \frac{x}{\log^2 x}, \tag{11}$$

where $c_3 > 0$ is a constant, see [RN, p. 265], [KV, p. 45], [MV, p. 91], [NT, p. 190]. Combining (10) and (11) puts the twin primes counting function in the form

$$a_1 \frac{x}{\log^2 x} \le \ \pi_2(x) \ \le b_1 \frac{x}{\log^2 x}, \tag{12}$$

where $a_1$, $b_1 > 0$ are constants. However, to establish the correct constant as conjectured

$$\pi_2(x) = 2c_2 \int_{1+\delta}^{x} \frac{dt}{\log^2 t} = 2c_2 li_2(x) + O(x^{1/2}), \tag{13}$$

where $\delta > 0$ is a small real number, and

$$c_2 = \prod_{p>2} \left(1 - (p-1)^{-2}\right) = .6601618158... \tag{14}$$

is the twin primes constant, much more efforts will be required. Perhaps, the probabilistic way of determining the constant is correct. The probabilistic argument produces asymptotic expression

$$\pi_2(x) = 2x \prod_{p>2} \left(1 - \frac{1}{(p-1)^2}\right) \prod_{2<p\le x^c} \left(1 - \frac{1}{p}\right)^2 = 2c_2 \frac{x}{\log^2 x} + O(\frac{x}{\log^3 x}), \tag{15}$$

where $c = e^{-\gamma} = e^{-.5772156649...} = 0.564307113...$ is a constant, and the number $1 - 1/p$ is the probability that the prime $p$ does not divides the primes $q$ nor $q + 2$. Confer [HW, p. 312], [KV, p. 38], and [ST, p. 180] for related details.

The general case of the Hardy-Littlewood conjecture has the prime pairs counting function

$$\pi_{2k}(x) = 2c_2 x \prod_{2<p\mid k} \left(\frac{p-1}{p-2}\right) + O(x^{1/2}), \tag{16}$$

where $k \ge 1$. The error term is conditional, see [KR, p. 7], and (13). The maximal order of magnitude is unconditionally bounded by

$$\pi_{2k}(x) \le c_{2k} \frac{x}{\log^2 x} \prod_{p\mid k} \left(1 + 1/p\right), \tag{17}$$

see [RN, p. 265]. The error term in this asymptotic formula is studied in [KT], and the numerical data for the change of signs $\pi_2(x) - 2c_2 li_2(x)$, here $li_2(x) = \int_2^x \log(t)^{-2}\, dt$ is the logarithm integral, is complied in [WA].

**Note 1.** This is reminiscent of the development of the prime number theorem: First, Chebyshev proved the inequality

$$a_0 \frac{x}{\log x} \le \pi(x) \le b_0 \frac{x}{\log x}, \tag{18}$$

where $a_0$, $b_0 > 0$ are constants, see [LA, p. 88], [LV, p. 105] for a proof. And about fifty years later, delaVallee Poussin and Hadamard proved the asymptotic formula

$$\pi(x) = x \prod_{2<p\le x^c} \left(1 - \frac{1}{p}\right) = \frac{x}{\log x} + O(\frac{x}{\log^2 x}), \tag{19}$$

where $c = e^{-\gamma} = e^{-.5772156649...} = 0.564307113...$, and the number $1 - 1/p$ is the probability that the prime $p$ does not divide the prime $q$.

The product representation of $\pi(x)$ is included here to demonstrate that the probabilistic argument produces the correct asymptotic formula as proved by analytical method. Perhaps, the only puzzling part is the constant $e^{-\gamma} = e^{-.5772156649...} = 0.564307113...$, consult [BW, p. 107], [HW] and [GV] for various levels of discussions.

### 2.3 Oscillations

A numerical experiment with the function $\pi_2(x) - 2c_2 li_2(x)$ has shown that it changes sign at much higher rate than $\pi(x) - li(x)$, see [WA]. The higher frequency of signs changes in $\pi_2(x) - 2c_2 li_2(x)$, as compared to $\pi(x) - li(x)$, appears to be due to the modulation effect of multiplying two density functions as in the case of twin primes. Specifically, the Gaussian density of primes is $\delta(x) = \psi(x)/x\log x$, and $li(x) = \int_2^x \delta(t)dt$. Similarly, the Gaussian density of twin primes is the product, (of sinusoidal functions),

$$\delta_T(x) = \frac{1}{\log x \log(x+2)} \left(1 - \sum_{|\rho|\le T} \frac{x^{\rho-1}}{\rho} + O(\log x / T)\right)\left(1 - \sum_{|\rho|\le T} \frac{(x+2)^{\rho-1}}{\rho} + O(\log x / T)\right), \tag{20}$$

where ρ are the nontrivial zeros of the zeta function $\zeta(s)$, $s \in \mathbb{C}$, on the critical strip $\{ s \in \mathbb{C} : 0 < \Re e(s) < 1 \}$, see [NW, p. 290], [TM, p. 177], and $li_2(x) = \int_2^x \delta_T(t)dt$. Equation (20) also seems to explain the smaller error term in $\pi_2(x) - 2c_2 li_2(x)$, by a factor of $\log x$, than in $\pi(x) - li(x)$.

### 2.4 Brun Constant

The Brun constant is defined by the twin primes series

$$B = \sum_p \left(\frac{1}{p} + \frac{1}{p+2}\right) = \left(\frac{1}{3} + \frac{1}{5}\right) + \left(\frac{1}{5} + \frac{1}{7}\right) + \left(\frac{1}{11} + \frac{1}{13}\right) + \cdots = 1.902160..., \tag{21}$$

this is a numerical approximation, more accurately, $1.830 < B < 2.347$, consult [KV, p. 12] for the advanced theory and computational techniques.

## 3. Germain Primes

A prime of the form $p = 2^a q + b$, where $q$ is prime, and $a$, $b \in \mathbb{Z}$, with $b$ odd are fixed integer parameters, is called a Germain prime. Usually, it has the form $p = 2q + 1$, see [RN, p. 331] for further discussion. Extremely large Germain primes are routinely computed by various workers in the field, see [PP], [WK], and the related literature. This work unfolds an argument in support of the existence of infinitely many Germain primes for any admissible fixed parameters $a$, $b \in \mathbb{Z}$. In particular, the cardinality of the subset of the usual Germain primes $\{ p, 2p + 1 : p \text{ is prime} \}$ is infinite.

### 3.1 The Main Result

***Theorem* 3.1.** For any fixed pair $a$, $b \in \mathbb{Z}$ of admissible integers the Diophantine equation $y = 2^a x + b$ has infinitely many prime solutions $x = p$, $y = q$.

Proof: Without loss in generality, fixed the pair $a = 1$, $b = 1$, and consider the weighted finite sum over the integers:

$$\begin{aligned}\sum_{n \le x} n\Lambda(n)\Lambda(2n+1) &= -\sum_{n \le x} n\Lambda(n) \sum_{d \mid 2n+1} \mu(d)\log d \\ &= -\sum_{d \le x} \mu(d)\log d \left( \sum_{n \le x,\ 2n \equiv -1 \bmod d} n\Lambda(n) \right),\end{aligned} \tag{22}$$

where $\Lambda(n) = -\sum_{d \mid n} \mu(d)\log d$. This follows from Lemma 6.1 in Section 6, and inverting the order of summation.

Other examples of inverting the order of summation are given in [MV, p. 35], [RE, p. 27], [RM, p. 216], [SP, p. 83], [TM, p. 36], et alii. The proof of Theorem 7.3 is quite similar to this, and it is instructive to review it.

Applying Lemmas 6.2, and 6.3 (with $a \equiv -1/2 \bmod q$, and $q = d$), and simplifying yield

$$\begin{aligned}\sum_{n \le x} n\Lambda(n)\Lambda(2n+1) &= -\sum_{d \le x,\ d \text{ odd}} \mu(d)\log d \left( \frac{1}{2\varphi(d)} x^2 + \Omega_{\pm}\left(x^{3/2} \log_3 x\right) \right) \\ &\ge -\sum_{d \le x,\ d \text{ odd}} \mu(d)\log d \left( \frac{1}{2d} x^2 \right),\end{aligned} \tag{23}$$

where $\log_3 x = \log\log\log x$ is the iterated logarithm, for infinitely many real numbers $x \ge 1$.

The restriction to odd $d = 2k + 1$ stems from the congruence relation $2n \equiv -1 \bmod d$ in (22). This restriction is equivalent to $\gcd(a, q) = 1$ with $a = 2k$, and $q = d$ in Lemmas 6.2, and 6.3. Here there is no restriction on the range of $d \le x$.

Applying Lemma 6.4 (for example, at $s = 1$, and $q_0 = 4$), the previous equation becomes

$$\begin{aligned}\sum_{n \le x} n\Lambda(n)\Lambda(2n+1) &\ge -c_1 x^2 \sum_{d \le x,\ d \text{ odd}} \frac{\mu(d)\log d}{d} \\ &= c_1 x^2 \left( c_2 + O\left( e^{-c(\log x)^{1/2}} \log x \right) \right) \\ &= c_3 x^2 + O\left( x^2 e^{-c(\log x)^{1/2}} \log x \right),\end{aligned} \tag{24}$$

where $c_1, c_2, c_3 > 0$ and $c > 0$ are nonnegative constants, for infinitely many real numbers $x \ge 1$.

Now, assume that there are finitely many Germain primes $p$, and $p + 2 < x_0$, where $x_0 > 0$ is a large constant. For example, $\Lambda(n)\Lambda(2n + 1) = 0$ for all $n > x_0$. Then

$$c_3 x^2 + O\left(x^2 e^{-c(\log x)^{1/2}} \log x\right) \le \sum_{n \le x} n\Lambda(n)\Lambda(2n+1) \le x_0^3 \tag{25}$$

for all real numbers $x \ge 1$. But this is a contradiction for infinitely many real numbers $x > x_0^2$. Ergo, there are infinitely many intervals $[1, x]$ that contain

$$\pi_G(x) \ge c_3 \frac{x}{\log^2 x} + O\left(xe^{-c(\log x)^{1/2}}\right) \tag{26}$$

Germain primes, see (28), as $x \ge 1$ tends to infinity. Quod erat demonstrandum. ■

### 3.2 Counting Function

The counting function of Germain primes is defined by

$$\pi_G(x) = \#\{\ p \le x : p, \text{ and } 2p+1 \text{ are primes}\ \}. \tag{27}$$

An application of partial summation to the finite sum yields the counting function

$$\begin{aligned}\pi_G(x) &= \sum_{n \le x} \frac{n\Lambda(n)\Lambda(2n+1)}{n\log(n)\log(2n+1)} - \sum_{n=p^k \text{ or } n+2=q^m \le x,\ k,m\ge 2} \frac{n\Lambda(n)\Lambda(2n+1)}{n\log(n)\log(2n+1)} \\ &= \int_2^x \frac{dR(t)}{t\log(t)\log(2t+1)} + O(x^{1/2}\log^2 x),\end{aligned} \tag{28}$$

where $R(x) = \sum_{n \le x} n\Lambda(n)\Lambda(2n+1)$. Evaluating the integral establishes a lower estimate of the Germain primes counting function

$$a\frac{x}{\log^2 x} \le \pi_G(x) \tag{29}$$

where $a > 0$ is a constant. This upper estimate as determined by sieve method, has the order of magnitude

$$\pi_G(x) \le b\frac{x}{\log^2 x}, \tag{30}$$

where $b > 0$ is a constant, see [RN, p. 265], [MV, p. 91], [NT, p. 190]. Combining (29) and (30) puts the Germain primes counting function in the form

$$a\frac{x}{\log^2 x} \le \pi_G(x) \le b\frac{x}{\log^2 x}, \tag{31}$$

where $a$, $b > 0$ are constants. The conjectured Germain primes counting function and its constant is the same as the twin primes constant, it est,

$$\pi_G(x) = 2c_2 \int_{1+\delta}^{x} \frac{dt}{\log^2 t} = 2c_2 li_2(x) + O(x^{1/2}), \tag{32}$$

where δ > 0 is a small real number, and

$$c_2 = \prod_{p>2} \left(1-(p-1)^{-2}\right) = .6601618158... \tag{33}$$

is the twin primes constant, much more efforts will be required. The error term is conditional, see [KR, p. 7].

### 3.3 Germain Constant

The same procedure used to prove the convergence of the series of the Brun constant is applied here to prove the convergence of the series of the new Germain constant. Let $q_n = 2p_n + 1$ be the $n$th Germain prime.

***Lemma* 3.2.** The series of the reciprocals of the $n$th Germain primes $p_n$ defined by the infinite sum

$$G = \sum_{p} \frac{1}{p} = \frac{1}{2} + \frac{1}{3} + \frac{1}{5} + \frac{1}{11} + \frac{1}{23} \cdots = 1.202203444... \tag{34}$$

converges.

Proof: Let $p_n$ be the $n$th Germain prime, and put $x = p_n$. Then, by (30), it follows that

$$n = \pi_G(p_n) \le c_1 p_n / \log^2 p_n, \tag{35}$$

where $c_1 > 0$ is a constant. Rearranging this inequality returns

$$\frac{1}{p_n} \le c_1 \frac{1}{n \log^2 p_n} \le c_1 \frac{1}{n \log^2 n}, \tag{36}$$

where $n \le p_n$. Therefore, the series

$$G = \sum_{n\ge1} \frac{1}{p_n} \le c_2 + c_1 \sum_{n\ge2} \frac{1}{n \log^2 n} < \infty, \tag{37}$$

where $c_2 \ge 0$ is a constant, converges, and defines the Germain constant. ■

The Brun constant is defined by the twin primes series

$$B = \sum_{p} \left(\frac{1}{p} + \frac{1}{p+2}\right) = \left(\frac{1}{3} + \frac{1}{5}\right) + \left(\frac{1}{5} + \frac{1}{7}\right) + \left(\frac{1}{11} + \frac{1}{13}\right) + \cdots = 1.902160..., \tag{38}$$

this is a numerical approximation, more accurately, $1.830 < B < 2.347$, consult [KV, p. 12] for the advanced theory and computational techniques.

Both the Brun constant, and the Germain constant are likely to be irrational numbers. However, it should be cautioned that a constant defined by an infinite sequence of prime numbers can be a rational number, exampli gratia, $\zeta(2)^2\zeta(4)^{-1} = \prod_{p\geq 2}\left((p^2+1)(p^2-1)^{-1}\right) = 5/2$.

### 3.3 Circle Method Comparison

Plenty of additive number theory problems have been solved by the circle method, [KM], [VN], [NT], [MT], and similar literature.

The circle method analysis of the Germain prime problem starts with a "weighted density function"

$$F(\alpha) = \left(\sum_{p\leq x}\log(p)e^{i2\pi p\alpha}\right)\left(\sum_{q<p}\log(q)e^{-i4\pi q\alpha}\right), \tag{39}$$

where $\Lambda(n)$ is the vonManglod function, see [MT, p. 32]. And the weighted counting function is

$$\begin{aligned}\sum_{p\leq x}\sum_{q<p}\log(q)\log(q) &= \int_0^1 F(\alpha)e^{i2\pi\alpha}d\alpha \\ &= \int_{\mathrm{M}} F(\alpha)e^{i2\pi\alpha}d\alpha + \int_m F(\alpha)e^{i2\pi\alpha}d\alpha \\ &= \Theta\frac{n}{2} + O\left(\frac{n}{\log n}\right) + \int_m F(\alpha)e^{i2\pi\alpha}d\alpha,\end{aligned} \tag{40}$$

where the limits of the integral are the major arcs *M*, and the minor arcs *m* of the unit interval [0, 1], and $\Theta(n)$ is the singular series. The exact value of the singular series is

$$\Theta(n) = \sum_{q\geq 1}\sum_{1\leq d\leq q}\frac{c_q(m)c_q(-2m)}{\varphi(q)^2}e^{-i2\pi d/q} = \prod_{p\geq 3}\left(1-\frac{1}{(p-1)^2}\right), \tag{41}$$

which is half the value of the twin prime constant, see [MT, 54].

The main term arising from the major arcs integral is known, but the minor arcs integral has not been estimated within acceptable order of magnitude, which is expected to be below the order of magnitude of the major arcs integral.

The local theory of Germain primes in short intervals and numerical data are given in and [WR], and [DB].

## 4. Primes In Quadratic Arithmetic Progressions

The cardinality of the subset of linear primes $\{ p = an + b : n \geq 1 \}$ defined by a linear polynomial $f(x) = ax + b \in \mathbb{Z}[x]$ over the integers, $\gcd(a, b) = 1$, was settled by Dirichlet as an infinite subset of prime numbers. The cardinality of the subset of quadratic primes $\{ p = an^2 + bn + c : n \geq 1 \}$ defined by certain irreducible quadratic polynomial $f(x) = ax^2 + bx + c \in \mathbb{Z}[x]$ over the integers, $\gcd(a, b, c) = 1$, is believed to be an infinite subset of prime numbers.

The quadratic primes conjecture claims that certain Diophantine equation $y = ax^2 + bx + c$ has infinitely many prime solutions $y = p$, as the integers $x = n \in \mathbb{Z}$ varies. More generally, the Bouniakowsky conjecture claims that for an irreducible $f(x) \in \mathbb{Z}[x]$ over the integers of fixed divisor $\text{div}(f) = 1$, and degree $\deg(f) \geq 2$, the Diophantine equation $y = f(x)$ has infinitely many prime solutions $y = p$ as the integers $x = n \in \mathbb{Z}$ varies.

The *fixed divisor* $\text{div}(f) = \gcd(f(\mathbb{Z}))$ of a polynomial $f(x) \in \mathbb{Z}[x]$ over the integers is the greatest common divisor of its image $f(\mathbb{Z}) = \{ f(n) : n \in \mathbb{Z} \}$ over the integers. The fixed divisor $\text{div}(f) = 1$ if $f(x) \equiv 0 \bmod p$ has $w(p) < p$ solutions for all primes $p < \deg(f)$, see [FI, p. 395].

Some irreducible polynomials have $\text{div}(f) \neq 1$, so these polynomials cannot generate infinitely many primes. For example, the irreducible polynomials $f_1(x) = x(x + 1) + 2$ has fixed divisor $\text{div}(f_1) = 2$, $f_2(x) = x(x + 1)(x + 2) + 3$ has fixed divisor $\text{div}(f_2) = 3$, and $f_3(x) = x^p - x + p$ has fixed divisor $\text{div}(f_3) = p$, so these polynomial cannot generate infinitely many primes. But, $f_4(x) = x^2 + 1$ has fixed divisor $\text{div}(f_3) = 1$, so it can generates infinitely many primes.

A reducible polynomial $f(x)$ of degree $d = \deg(f)$ can generates up to $d \geq 0$ primes. For example, the reducible polynomial $f(x) = (g(x) \pm 1)(x^2 + 1)$ can generates up to $\deg(g) \geq 0$ primes for any $0 \neq g(x) \in \mathbb{Z}[x]$ of degree $\deg(g) \geq 0$.

Detailed discussions appear in [RN, p. 387], [GL, p. 17], [FI, p. 395], [NW, p. 405], [PZ, p. 33], [WK], [WS], [PP], and related topics in [BR], [MA], et alii.

This work provides the details of a possible proof of this quadratic primes counting problem. In particular, the cardinality of the simplest subset of quadratic primes $\{ p = n^2 + 1 : n \geq 1 \}$ is infinite.

### 4.1 The Main Result

***Theorem* 1.** Let $f(x) = ax^2 + bx + c \in \mathbb{Z}[x]$ be an irreducible polynomial over the integers of fixed divisor $\text{div}(f) = 1$. Then, the Diophantine equation $p = an^2 + bn + c$ has infinitely many primes $p$ solutions as $n \in \mathbb{Z}$ varies over the integers.

Proof: Without loss in generality, let $a = 1$, $b = 0$, $c = 1$, and consider the weighted finite sum over the integers:

$$\begin{aligned} \sum_{n^2 \leq x} n^2 \Lambda(n^2 + 1) &= -\sum_{n^2 \leq x} n^2 \sum_{d \mid n^2+1} \mu(d) \log d \\ &= -\sum_{d \leq x+1} \mu(d) \log d \left( \sum_{n \leq x+1,\ n^2 \equiv -1 \bmod d} n^2 \right), \end{aligned} \tag{42}$$

where $\Lambda(n) = -\sum_{d \mid n} \mu(d) \log d$ . This follows from Lemma 6.1 in Section 6, and inverting the summation. Other examples of inverting the order of summation are given in [MV, p. 35], [RE, p. 27], [RM, p. 216], [SP, p. 83], [TM, p. 36], et alii.

Applying Lemma 8.1, with ($a = -1$ and $q = d$), and simplifying yield

$$\begin{aligned}\sum_{n^2 \le x} n^2\Lambda(n^2+1) &\ge -\sum_{d \le x+1,\ n^2 \equiv -1 \bmod d} \mu(d)\log d\left(\frac{2^W}{3d^2}\left(x^{1/2}+1\right)^3 + O\left(\frac{1}{d}(x^{1/2}+1)^2\right)\right) \\ &\ge -c_1x^{3/2}\sum_{d \le x+1,\ p \mid d \Rightarrow p \equiv 1 \bmod 4} \frac{\mu(d)\log d}{d^2} + O\left(x\left|\sum_{d \le x+1,\ p \mid d \Rightarrow p \equiv 1 \bmod 4} \frac{\mu(d)\log d}{d}\right|\right),\end{aligned} \tag{43}$$

where $c_1 > 0$ is a constant, and $2^W \ge 2$. The congruence $n^2 \equiv -1 \bmod d$ is replaced with the equivalent restriction $p \mid d \Rightarrow p \equiv 1 \bmod 4$.

Applying Lemma 8.2, (for example, at $s = 2$, for $q_0 = 4$, and $s = 1$, for $q_0 = 4$), the previous equation becomes

$$\begin{aligned}\sum_{n^2 \le x} n^2\Lambda(n^2+1) &\ge c_1x^{3/2}\left(c_2 + O\left(x^{-1}e^{-c(\log x)^{1/2}}\log x\right)\right) + O\left(xO\left(c_3 + e^{-c(\log x)^{1/2}}\log x\right)\right) \\ &= c_4x^{3/2} + O(x),\end{aligned} \tag{44}$$

where $c_1$, $c_2$, $c_3$, $c_4 > 0$, and $c > 0$ are nonnegative constants. Now, assume that there are finitely many quadratic primes $p = n^2 + 1 < x_0$, where $x_0 > 0$ is a large constant. For example, $\Lambda(n^2 + 1) = 0$ for all $n > x_0$. Then

$$c_3x^{3/2} + O(x) \le \sum_{n^2 \le x} n^2\Lambda(n^2+1) \le x_0^2 \tag{45}$$

for any real number $x \ge 1$. But this is a contradiction for all sufficiently large real number $x > x_0^2$. Lux mentis lux orbis. ■

Some irreducible polynomials $f(x) = ax^2 + bx + c \in \mathbb{Z}[x]$ of fixed divisor dev($f$) = 1 can be converted to an equivalent case of the form $g(x) = x^2 + d \in \mathbb{Z}[x]$ by means of algebraic or analytical manipulations. The equivalent problem is then handled as the case $x^2 + 1 \in \mathbb{Z}[x]$ presented above. The general case $f(x) = ax^2 + bx + c \in \mathbb{Z}[x]$ can be handled via the weighted finite sum

$$\sum_{n \le x}(an^2+bn)\Lambda(an^2+bn+c) = a\sum_{n \le x} n^2\Lambda(an^2+bn+c) + b\sum_{n \le x} n\Lambda(an^2+bn+c), \tag{46}$$

or other similar finite sums. It should also be noted that other weighted finite sums over the integers such as

$$\sum_{n \le x}\varphi(n)\Lambda(n^2+1), \qquad \sum_{n \le x}\sigma(n)\Lambda(n^2+1), \qquad \sum_{n \le x}\tau(n)\Lambda(n^2+1), \tag{47}$$

where φ, σ, and τ are the totient function, and the divisors functions respectively, yield precisely the same result. The analysis of these later finite sums are quite similar to the Titchmarsh divisor problem over quadratic arithmetic progressions, see [CM, p. 172] for an introduction to this problem.

### 4.2. Counting Functions

The expected distribution of the quadratic primes has its systematic development in the middle of the last century, and it is described in the well-known series of lectures Partitio Numerorum of Hardy and Littlewood. This claim will be referred to as the *quadratic primes conjecture*. This conjecture and other generalizations are discussed in [RN, p. 406], [GL, p. 25], [NW, p. 342], and similar sources. Related works are given in [BR], [MA], [GM], and [JW].

***Quadratic Primes Conjecture.*** Let $f(x) = ax^2 + bx + c \in \mathbb{Z}[x]$ be an irreducible polynomial over the integers, and assume that

(i) $\gcd(a, b, c) = 1$, (ii) $b^2 - 4ac$ is not a square in $\mathbb{Z}$, (iii) $\gcd(a + b, c) = 2d + 1$ is odd.

Then, for a sufficiently large real number $x \geq 1$, the number of primes of the form $p = an^2 + bn + c \leq x$ has the asymptotic formula

$$\pi_f(x) = \frac{\varepsilon C_f}{\sqrt{a}} \frac{x^{1/2}}{\log x} \prod_{p \mid \gcd(a,b)} \left(1 + (p-1)^{-1}\right) + O\left(\frac{x^{1/2}}{\log^2 x}\right), \tag{48}$$

where the constants are

$$C_f = \prod_{p>2,\ \gcd(a,p)=1} \left(1 + \chi(b^2 - 4ac)(p-1)^{-1}\right), \quad \text{and} \quad \varepsilon = \begin{cases} 1 & \text{if } a + b \equiv 1 \bmod 2, \\ 2 & \text{if } a + b \equiv 0 \bmod 2, \end{cases} \tag{49}$$

here $\chi(t)$ denotes the quadratic symbol modulo $p$.

Conditions (i) and (ii) imply that $f(x)$ is irreducible, and condition (iii) implies that the fixed divisor $\text{div}(f) = \gcd(f(\mathbb{Z})) = 1$, [FI, p. 395].

The simplest case $f(x) = x^2 + 1 \in \mathbb{Z}[x]$ reduces to

$$\pi_f(x) = \frac{x^{1/2}}{\log x} \prod_{p>2} \left(1 + (-1)^{(p-1)/2}(p-1)^{-1}\right) + O\left(\frac{x^{1/2}}{\log^2 x}\right) \tag{50}$$

for a sufficiently large real number $x \geq 1$. The constant $C_f = 1.37281346\ldots$, [CW, p. 13].

The complete proof, including the verification of the constant $C = \varepsilon a^{-1/2} C_f \prod_{p \mid \gcd(a,b)} (1 + (p-1)^{-1})$, of any of these cases is not expected to be completed any time soon. However, a lower estimate of the correct order of magnitude can be derived by from Theorem 1 by partial summation. Specifically, for the primes $p = n^2 + 1 \leq x$, the quadratic primes counting function is given by

$$\begin{aligned}\pi_f(x) &= \sum_{n \le x} \frac{n^2\Lambda(n^2+1)}{n^2\log(n^2+1)} - \sum_{n=p^k \le x^{1/2},\ k\ge 2} \frac{n^2\Lambda(n^2+1)}{n^2\log(n^2+1)} \\ &= \int_2^{x^{1/2}} \frac{dR(t)}{t^2\log(t^2+1)} + O(x^{1/4}\log^2 x),\end{aligned} \tag{51}$$

where $R(x) = \sum_{n^2 \le x} n^2\Lambda(n^2+1)$. This leads to a lower estimate of the quadratic primes counting function

$$\pi_f(x) = \#\left\{ p = n^2+1 \le x : n \ge 1 \right\} \ge a\frac{x^{1/2}}{\log^2 x} + O(\frac{x^{1/2}}{\log^2 x}), \tag{52}$$

with $a > 0$ constant. The upper estimate

$$\pi_f(x) \le b\frac{x^{1/2}}{\log^2 x}, \tag{53}$$

was proved in [GM], this also follows from (1) and the Brun-Titchmarsh theorem. Combining (52) and (53) puts the quadratic primes counting function in the form

$$a\frac{x^{1/2}}{\log^2 x} \le \pi_f(x) \le b\frac{x^{1/2}}{\log^2 x}, \tag{54}$$

where $a, b > 0$ are constants.

### 4.3 Optimization Problem

The problem of finding quadratic polynomials with high densities of prime numbers in very short interval $[0, x]$ has intrigued people for a long time. The Euler polynomial $f(x) = x^2 + x + 41$ has the highest prime density known for small $x < 41$. This optimization problem, which seeks the best constant $C = \varepsilon a^{-1/2}C_f\prod_{p \mid \gcd(a,b)}(1+(p-1)^{-1})$ in equation (48) for some in $f(x) \in \mathbb{Z}[x]$, is studied in [JW], [GM] and many other.

### 4.4 Quadratic Twin Primes

The quadratic twin primes conjecture calls for a subset of infinitely many pairs of the forms $p = n^2 + 1$, and $p + 2 = n^2 + 3$, see [GL, p. 27], [RN], [NW, p. 343].

***Quadratic Twin Primes Conjecture.*** There are infinitely many quadratic twin primes $n^2 + 1$, and $n^2 + 3$. The counting function has the asymptotic formula

$$\pi_{QTP}(x) = 6\frac{x^{1/2}}{\log x}\prod_{p\ge 5}\frac{p(p-w(p))}{(p-1)^2} + O(\frac{x^{1/2}}{\log^2 x}), \tag{55}$$

where $w(p)$ denotes the number of quadratic residues mod $p$ in the set $\{ -1, -3 \}$.

It seems plausible that the correct order of magnitude for this problem can be handled in the same way as the linear twin primes problem. Id est, estimating the weighted sum

$$\sum_{n^2 \le x} n^2 \Lambda(n^2+1)\Lambda(n^2+3) \tag{56}$$

by combining the previous analysis of the linear twin primes in Theorem 2.1, and the analysis done here for the quadratic primes in Theorem 3.1.

## 5. Primes In Linear Arithmetic Progressions

The subset of primes in linear arithmetic progressions $\{ p = an + b : n \ge 1 \}$, $\gcd(a, b) = 1$, was studied by Dirichlet using complex analysis. Later, Selberg used elementary methods to prove the same result, see [SB].

The cardinality of the subset of primes in many linear arithmetic progressions can be determined by Euclidean type proofs, for examples, $a = 4$, 6 and so on, but there many cases that cannot be determined by Euclidean type proofs, for examples, $\{ p = 5n + 3 : n \ge 1 \}$. An earlier result of Schur shows that a Euclidean type proof exists if $b^2 \equiv 1 \bmod a$. The converse was recently proved, see [MR].

### 5.1 The Main Result

A simpler elementary and characterfree proof of the existence of infinitely many primes in an arithmetic progression $\{ p = qn + a : n \ge 1 \}$ is sketched below. However, the constant of the main term remains undetermined.

***Theorem* 5.1.** Let $a \ge 1$, and $q \ge 1$ be fixed integers, $\gcd(a, q) = 1$. Then, the Diophantine equation $y = qx + a$ has infinitely many prime solutions $y = p$ as $x \in \mathbb{Z}$ varies over the integers.

Proof: Consider the weighted finite sum over the integers:

$$\begin{aligned} \sum_{n \le x} n\Lambda(qn+a) &= -\sum_{n \le x} n \sum_{d \mid qn+a} \mu(d)\log d \\ &= -\sum_{d \le x} \mu(d)\log d \left( \sum_{n \le x,\ n \equiv -aq^{-1} \bmod d} n \right), \end{aligned} \tag{57}$$

where $\Lambda(n) = -\sum_{d \mid n} \mu(d)\log d$. This follows from Lemma 6.1 in Section 6, and inverting the order of summation.

Applying Lemma 8.3, with ($a_0 = -a/q$ and $q_0 = d$), and simplifying yield

$$\begin{aligned} \sum_{n \le x} n\Lambda(qn+a) &= -\sum_{d \le x,\ \gcd(d,q)=1} \mu(d)\log d \left( \frac{1}{2d} x^2 + O(x) \right) \\ &= -c_1 x^2 \sum_{d \le x,\ \gcd(q,d)=1} \frac{\mu(d)\log d}{d} + O\left( x \left| \sum_{d \le x,\ \gcd(d,q)=1} \mu(d)\log d \right| \right), \end{aligned} \tag{58}$$

where $c_1 > 0$ is a constant. The restriction $n \equiv -a/q \bmod d$ is replaced with the equivalent restriction $\gcd(d, q) = 1$.

Applying Lemma 6.4 (for example, at $s = 1$ for $q = q$), and Lemma 6.5, the previous equation becomes

$$\sum_{n \le x} n\Lambda(qn+a) = c_1 x^2 \left(c_2 + O\left(e^{-c(\log x)^{1/2}} \log x\right)\right) + O\left(xO\left(xe^{-c(\log x)^{1/2}} \log x\right)\right)$$
$$= c_3 x^2 + O\left(x^2 e^{-c(\log x)^{1/2}} \log x\right), \tag{59}$$

where $c_1$, $c_2$, $c_3 > 0$ and $c > 0$ are nonnegative constants. Further, by partial summation, it follows that

$$\pi(x,a,q) = \sum_{n \le x} \frac{n\Lambda(qn+a)}{n\log(qn+a)} + O(x^{1/2}) \ge c_3 \frac{x}{\log x} + O\left(xe^{-c(\log x)^{1/2}}\right) \tag{60}$$

for all large real number $x \ge 1$. ■

As shown above, the determination of the constant $c_3 = 1/\varphi(q)$ is not required to confirm that the cardinality of the subset of primes in a linear arithmetic progression $\{ p = an + b : n \ge 1 \}$, $\gcd(a, b) = 1$, is infinite. This proof is based on the same foundation as the elementary proof of the prime number theorem on arithmetic progression, confer [SL]. However, it is somewhat simpler.

## 6. Elementary Foundation

The elementary underpinnings of Theorems 2.1, 3.1, 4.1 and 5.1 are assembled here. The basic definitions of several number theoretical functions, and a handful of Lemmas are recorded here.

### 6.1 Formulae for the vonMangoldt Function

Let $n \in \mathbb{N} = \{ 0, 1, 2, 3, \ldots \}$ be an integer. The Mobius function is defined by

$$\mu(n) = \begin{cases} (-1)^v & \text{if } n = p_1 p_2 \cdots p_v, \\ 0 & \text{if } n \ne \text{squarefree}, \end{cases} \tag{61}$$

where the $p_i$ are primes. And the vonMangoldt function, and its more general version are defined by

$$\Lambda(n) = \begin{cases} \log p & \text{if } n = p^m, m \ge 1, \\ 0 & \text{if } n \ne p^m, \end{cases} \quad \text{and} \quad \Lambda_k(n) = \begin{cases} \log^k p & \text{if } \omega(n) \le k, \\ 0 & \text{if } \omega(n) > k, \end{cases} \tag{62}$$

for a fixed integer $k \ge 1$. A few of the various representations of the vonMangoldt function (62) are shown here. There is a simple proof of the identity (63) via the Mobius inversion formula, but an analytical proof is given to show its connection to the zeta function and beyond.

***Lemma* 6.1.** Let $k \ge 1$ be a fixed integer. Let $n \ge 1$ be an integer, and let $\Lambda_k(n)$ be the vonMangoldt function. Then

$$\text{(i)} \quad \Lambda(n) = -\sum_{d \mid n} \mu(d) \log d, \qquad \text{(ii)} \quad \Lambda_k(n) = \sum_{d \mid n} \mu(d) \log(n/d)^k. \tag{63}$$

Proof: The logarithm derivative of the zeta function satisfies

$$-\frac{\zeta'(s)}{\zeta(s)}=\sum_{n\geq 1}\frac{\Lambda(n)}{n^s}. \tag{64}$$

On the other hand, the convolution of the Dirichlet series of the functions $1/\zeta(s)$, and $-\zeta'(s)$ respectively satisfies

$$-\frac{\zeta'(s)}{\zeta(s)}=\left(\sum_{n\geq 1}\frac{\log n}{n^s}\right)\left(\sum_{n\geq 1}\frac{\mu(n)}{n^s}\right)=\sum_{n\geq 1}\frac{\sum_{d|n}\mu(d)\log n/d}{n^s}. \tag{65}$$

Lastly, matching the coefficients of the powers of $n^{-s}$ in (64) and (65) confirms the claim. ■

Various other identities and approximations of the vonMangoldt function are used in the literature. These include short formulas such as

$$\Lambda_R(n)=\sum_{d\,|\,n,d\leq R}\mu(d)\log R/d\ , \tag{66}$$

see [GP], the Vaughan identity, see [HN, p. 25], and the harmonic expansion

$$\Lambda(n)=-\sum_{q\geq 1}\frac{c_n(q)}{q}, \tag{67}$$

where $c_n(q)=\sum_{\gcd(x,n)=1}e^{i2\pi qx/n}$, see [MV, p. 114]. This later identity leads to the conjectured twin prime formula

$$\sum_{n\leq x}\Lambda(n)\Lambda(n+2)=-\sum_{q\geq 1}\frac{\sum_{n\leq x}\Lambda(n)c_{n+2}(q)}{q}=2c_2x+o(x), \tag{68}$$

where $c_2$ is the twin prime constant, see (15), and large real numbers $x \geq 1$.

Extensive details of other identities and approximations of the vonMangoldt are discussed in [FI], et alii.

### 6.2 Weighted Finite Sum Over Arithmetic Progressions

This result is an important part of Theorems 2.1, and 3.1. The proof uses elementary methods. Other analytical methods are possible and provide alternative proofs of the twin primes conjecture, and related primes counting problems.

The weighted finite sum of interest has the explicit formula

$$\sum_{n\leq x}n\Lambda(n)=\frac{x^2}{2}+\sum_{\rho}\frac{x^{\rho+1}}{\rho(\rho+1)}+\frac{\zeta'(0)}{\zeta(0)}-\frac{\zeta'(1)}{\zeta(1)}-\sum_{n\geq 1}\frac{x^{1-2n}}{2n(2n-1)}\ , \tag{69}$$

where ρ are the nontrivial zeros of the zeta function $\zeta(s)$, $s \in \mathbb{C}$, on the critical strip $\{ s \in \mathbb{C} : 0 < \Re e(s) < 1 \}$, see [EL, p. 162], [MV, p. 420]. This formula immediately shows the oscillatory nature of this finite sum, and implies that $\sum_{n \le x,\ n \equiv a \bmod q} n\Lambda(n) \ge x^2 / 2\varphi(q)$ infinitely often as the real number $x \ge 1$ tends to infinity. This is proven below.

***Lemma* 6.2.** Let $a \ge 1$, and $q \ge 1$ be integers such that and $\gcd(a, q) = 1$. Let $x \ge 1$ be a large real number. Then

(i) $$\sum_{n \le x,\ n \equiv a \bmod q} n\Lambda(n) = \frac{1}{2\varphi(q)} x^2 + O(x^2 e^{-c(\log x)^{1/2}}) , \tag{70}$$

where $q = O((\log x)^B)$, $B > 0$, and $c > 0$ constants.

(ii) $$\sum_{n \le x,\ n \equiv a \bmod q} n\Lambda(n) = \frac{1}{2\varphi(q)} x^2 + \Omega_{\pm}(x^{3/2} \log\log\log x),$$

where $q = O( x^{1/2})$.

(iii) $$\sum_{n \le x,\ n \equiv a \bmod q} n\Lambda(n) = \frac{1}{2\varphi(q)} x^2 + O(x^{3/2}) ,$$

where $q = O(x^{1/2})$, and assuming the Riemann hypothesis.

Proof (ii): By Theorems 7.2, and 7.3 the Chebychev psi function satisfies

$$\psi(x,a,q) = \frac{x}{\varphi(q)} + \Omega_{\pm}(x^{1/2} \log\log\log x) \tag{71}$$

on the arithmetic progression $\{ qn + a : n \ge 1, \text{ and } \gcd(a, q) = 1 \}$. Now, use the measure $d\psi(t,a,q)$ to evaluate the integral representation

$$\begin{aligned} \sum_{n \le x,\ n \equiv a \bmod q} n\Lambda(n) &= \int_1^x t\, d\psi(t,a,q) \\ &= \frac{1}{\varphi(q)} x^2 - \int_1^x \psi(t,a,q)\,dt \\ &= \frac{1}{2\varphi(q)} x^2 + \Omega_{\pm}\left(x^{3/2} \log\log\log x\right) \end{aligned} \tag{72}$$

for all large real numbers $x \ge 1$. The proofs of (i) and (iii) are similar. ■

A different proof of (iii) appears in [NW, p. 330].

The *level of distribution* of the modulo $q$ is the maximal value $q \le x^{\theta}$ admissible for some $0 \le \theta \le 1$. The Elliot-

Haberstram conjecture states that θ = 1 – ε with ε > 0 arbitrarily small. Explanation of the concept of level of distribution appears in [SD], [FI, 405], the advanced theory appears in [FR], and similar literature.

***Lemma* 6.3.** Let $a \geq 1$, and $q \geq 1$ be integers such that and $\gcd(a, q) = 1$. Let $x \geq 1$ be a large real number. Then

$$\sum_{n \leq x,\ n \equiv a \bmod q} n\Lambda(n) \geq \frac{1}{2q} x^2, \tag{73}$$

where $q \leq x$, for infinitely many real numbers $x \geq 1$.

Proof: By Lemma 6.2, the omega result

$$\begin{aligned} \sum_{n \leq x, n \equiv a \bmod q} n\Lambda(n) &= \frac{1}{2\varphi(q)} x^2 + \Omega_+(x^{3/2} \log_3 x) \\ &\geq \frac{1}{2\varphi(q)} x^2 + c_0 x^{3/2} \log_3 x > 0, \end{aligned} \tag{74}$$

where $c_0 > 0$ is a constant, and $\log_3 x = \log\log\log x$ is the iterated logarithm, holds for all $d \leq x$, and infinitely many real numbers $x \geq 1$. Here, no information on the level of distribution of the modulo $d$ is required in this application. In view of these data, the last equation becomes

$$\begin{aligned} \sum_{n \leq x} n\Lambda(n) &\geq \frac{1}{2\varphi(q)} x^2 + c_0 x^{3/2} \log_3 x \\ &\geq \frac{1}{2\varphi(q)} x^2 \\ &\geq \frac{1}{2q} x^2 \end{aligned} \tag{75}$$

for infinitely many real numbers $x \geq 1$. ■

This estimate will be utilized in Theorems 2.1, and 3.1.

### 6.3 Twisted Finite Sums

***Lemma* 6.4.** Let $s \geq 1$ be an integer, let $\mu(n)$ be the Mobius function, and let $x \geq 1$ be a sufficiently large real number. Then

$$\text{(i)} \quad \sum_{n \leq x} \frac{\mu(n) \log n}{n^s} = \frac{\zeta'(s)}{\zeta(s)^2} + O\left(x^{1-s} e^{-c(\log x)^{1/2}} \log x\right), \tag{76}$$

$$\text{(ii)} \quad \sum_{n \leq x,\ n \equiv a \bmod q} \frac{\mu(n) \log n}{n^s} = \frac{L'(s, \chi_0)}{L(s, \chi_0)^2} + O\left(x^{1-s} e^{-c(\log x)^{1/2}} \log x\right),$$

(iii) $$\sum_{n \le x,\ \gcd(a,q)=1} \frac{\mu(n)\log n}{n^s} = -C_0 + O\left(x^{1-s} e^{-c(\log x)^{1/2}} \log x\right),$$

where $a \geq 1$, $q \geq 1$ is a pair of fixed integers, gcd($a$, $q$) = 1, $C_0 > 0$, $c > 0$ are constants, and $\chi_0$ is the principal character modulo $q$.

Proof (i): Write the finite sum in the form

$$\sum_{n \le x} \frac{\mu(n)\log n}{n^s} = \sum_{n \ge 1} \frac{\mu(n)\log n}{n^s} - \sum_{n > x} \frac{\mu(n)\log n}{n^s}. \tag{77}$$

The constant is expressed in terms of the logarithmic derivative of the zeta function, id est,

$$\frac{d}{ds}\frac{1}{\zeta(s)} = -\sum_{n \ge 1} \frac{\mu(n)\log n}{n^s} = -\frac{\zeta'(s)}{\zeta(s)^2}. \tag{78}$$

Now, since the zeta function $\zeta(s) = \sum_{n \ge 1} n^{-s}$ is a decreasing function on the real half line $\Re e(s) > 1$, the derivative $\zeta'(s) = -\sum_{n \ge 1} (\log n) n^{-s}$ is negative on the real half line $s > 1$. The result follows from these data. In (ii) repeat the routine for the L-function $L(s, \chi_0)$, the constant $C_0 = q/\varphi(q)$ appears in [MV, p. 185]. ■

The zeta function, and its derivative have the power series expansions

$$\zeta(s) = \frac{1}{s-1} + \gamma_0 - \gamma_1(s-1) + \frac{1}{2}\gamma_2(s-1)^2 - \cdots, \tag{79}$$

$$\zeta'(s) = -\frac{1}{(s-1)^2} - \gamma_1 + \gamma_2(s-1) - \frac{1}{2}\gamma_3(s-1)^2 + \cdots$$

where $\gamma_k$ is the $k$th Stieltjes constant, see [DL, 25.2.4]. These series can be used to compute numerical approximations of the constants $\kappa_s = -\zeta'(s)/\zeta(s)^2$ for $s \geq 1$. Exampli gratia, the logarithmic derivative at $s = 1$ is

$$\frac{\zeta'(1)}{\zeta(1)^2} = \lim_{s \to 1} \frac{(s-1)^2 \zeta'(s)}{(s-1)^2 \zeta(s)^2} = -1, \tag{80}$$

and

$$\sum_{n \le x} \frac{\mu(n)\log n}{n} = \frac{\zeta'(1)}{\zeta(1)^2} + O\left(e^{-c(\log x)^{1/2}} \log x\right) = -1 + O\left(e^{-c(\log x)^{1/2}} \log x\right), \tag{81}$$

where $c > 0$ is a constant. Similarly, the logarithmic derivative at $s = 2$ is exactly

$$\frac{\zeta'(2)}{\zeta(2)^2} = -\frac{(\pi^2/6)(\gamma + \log 2\pi - 12\log C)}{\pi^4/36} = -\frac{6(\gamma + \log 2\pi - 12\log C)}{\pi^2}, \tag{82}$$

where $\gamma$ is Euler constant, and $C$ is Glaisher constant, see [DL, 5.17.7]. The combined result is

$$\sum_{n\leq x}\frac{\mu(n)\log n}{n^2}=\frac{\zeta'(2)}{\zeta(2)^2}+O\left(e^{-c(\log x)^{1/2}}\log x\right)=-\frac{6(\gamma+\log 2\pi-12\log C)}{\pi^2}+O\left(e^{-c(\log x)^{1/2}}\log x\right). \tag{83}$$

In the other case (ii), the L-function and its derivative are defined by

$$L(s,\chi)=\sum_{n\geq 1}\frac{\chi(n)}{n^s}, \quad \text{and} \quad L'(s,\chi)=-\sum_{n\geq 1}\frac{\chi(n)\log n}{n^s}, \tag{84}$$

where $\chi$ is a character modulo $q$, see [DL, 5.15.1].

***Lemma* 6.5.** Let $x \geq 1$ be a sufficiently large real number, and let $\mu(n)$ be the Mobius function. Then

$$\text{(i)}\ \sum_{n\leq x}\mu(n)=O\left(xe^{-c(\log x)^{1/2}}\right), \qquad \text{(ii)}\ \sum_{n\leq x}\mu(n)\log n=O\left(xe^{-c(\log x)^{1/2}}\log x\right), \tag{85}$$

$$\text{(i)}\ \sum_{n\leq x}\mu(n)=\Omega_{\pm}\left(x^{1/2}\right),$$

where $c > 0$ is a constant.

Proof (i) : Sketches of the proof are given in [MV, p. 182], and [RM, p. 316]. To prove (ii), evaluate the Stieltjes integral

$$\sum_{n\leq x}\mu(n)\log n=\int_1^x(\log t)\,dM(t), \tag{86}$$

where $M(x)=\sum_{n\leq x}\mu(n)=O(xe^{-c(\log x)^{1/2}})$. ■

## 7. Prime Numbers Theorems

The primes counting function is defined by

$$\pi(x)=\#\{\ p\leq x : p\ \text{ is prime}\ \}, \tag{87}$$

and the primes counting function on an the arithmetic progression $\{ qn + a : n \geq 1\}$, $\gcd(a, q) = 1$, is defined by

$$\pi(x,a,q)=\#\{\ p\equiv a \bmod q\leq x : p\ \text{ is prime}\ \}. \tag{88}$$

The Chebychev psi function is defined by

$$\psi(x)=\sum_{n\leq x}\Lambda(n) \quad \text{or} \quad \psi(x,a,q)=\sum_{n\leq x, n\equiv a \bmod q}\Lambda(n) \tag{89}$$

over the integers $\mathbb{N} = \{ 0, 1, 2, 3, \ldots \}$, or over the arithmetic progression $q\mathbb{N} + a = \{ qn + a : n \geq 0 \}$ respectively.

The earliest form of the prime number theorem was Chebychev inequality with the correct order of magnitude. The original proof is much longer than the current versions.

***Lemma* 7.1.** There is a pair of constants $a < b$ such that

$$a\frac{x}{\log x} \le \pi(x) \le b\frac{x}{\log x}. \tag{90}$$

The proof is available in [LA, p. 88], [LV, p. 105], and almost any textbook in Number Theory.

Around fifty years after the prime number theorem was proved by means of complex variable, two elementary proofs were discovered, one by Erdos, and another by Selberg. These proofs are based on the Selberg formula

$$\sum_{p \le x} \log^2 p + \sum_{pq \le x} \log p \log q = 2x\log x + O(x). \tag{91}$$

The elementary proofs of the prime number theorem are given in [ES], [DD], [SB], and other authors. Another completely elementary proof of the prime number theorem is provided below.

***Theorem* 7.2.** Let $x \ge 1$ be a large real number. Then

$$\pi(x) = \frac{x}{\log x} + O(\frac{x}{\log^2 x}). \tag{92}$$

Proof: Let $x \ge 1$ be a large integer, and write

$$x! = \prod_{p \mid x} p^{v_p}, \tag{93}$$

where $v_p = \sum_{n \ge 1}\left[x/p^n\right]$. Taking logarithm yields

$$\log x! = \sum_{p \le x}\left(\left[\frac{x}{p}\right] + \left[\frac{x}{p^2}\right] + \left[\frac{x}{p^3}\right] + \cdots\right)\log p = x\sum_{p \le x}\frac{\log p}{p} + O\left(x\right), \tag{94}$$

see [LV, p. 102]. Now, the derivation of an estimate of the finite sum on the right side requires no knowledge of prime numbers: An integral approximation yields

$$\log x! = \sum_{n \le x} \log n = x\log x - x + C_0 + O(\log x) \tag{95}$$

where $C_0 > 0$ is a constant, or invoke the Sterling formula

$$\log\Gamma(z) = \left(z - 1/2\right)\log z - z - \log(2\pi)/2 + O(1/|z|), \tag{96}$$

where $\Gamma(z+1) = z\Gamma(z)$ is the gamma function, and $z$ is a complex number. Comparing (94) and (95) yields

$$\sum_{p\le x}\frac{\log p}{p}=\log x-1+C_1+O\left(\frac{\log x}{x}\right), \tag{97}$$

where $C_1 > 0$ is a constant. By partial summation, the number of primes is given by

$$\begin{aligned}\pi(x)&=\sum_{n\le x}\frac{\Lambda(n)}{n}\frac{n}{\log n}-\sum_{n=p^k\le x,\ k\ge2}\frac{n\Lambda(n)}{n\log n}\\&=\int_2^{x^{1/2}}\frac{t}{\log t}dR(t)+O(x^{1/2}\log x)\\&=\frac{x}{\log x}+O(\frac{x}{\log^2 x}),\end{aligned} \tag{98}$$

where $R(x)=\sum_{n\le x}\log p/p$. ■

***Theorem* 7.3.** Let $x \ge 1$ be a large real number. Then

(i) $\pi(x)=li(x)+O\left(xe^{-c(\log x)^{1/2}}\right)$, (99)

(ii) $\psi(x)=x+O\left(xe^{-c(\log x)^{1/2}}\right)$,

where $c > 0$ is an absolute nonnegative constant.

Proof (ii): Let (( $x$ )) be the fractional part function. By Lemma 6.1, the psi function $\psi$ can be rewritten as

$$\begin{aligned}\sum_{n\le x}\Lambda(n)&=-\sum_{n\le x}\sum_{d\,|\,n}\mu(d)\log d\\&=-\sum_{d\le x}\mu(d)\log d\left(\sum_{n\le x/d}1\right)\\&=-\sum_{d\le x}\mu(d)\log d\left(\frac{x}{d}-((\ x/d\ ))\right),\end{aligned} \tag{100}$$

where (( $x$ )) is the exact error incurred in the approximation of the finite sum $\sum_{n\le x/d}1\ge0$. Applying Lemmas 6.4, and 6.5 and simplifications return

$$\begin{aligned}\sum_{n\le x}\Lambda(n)&=-x\sum_{d\le x}\frac{\mu(d)\log d}{d}+O\left(\left|\sum_{d\le x}\mu(d)\log d\big((x/d)\big)\right|\right)\\&=x\left(1+O\left(e^{-c(\log x)^{1/2}}\log x\right)\right)+O\left(xe^{-c(\log x)^{1/2}}\log x\right)\\&=x+O\left(xe^{-c(\log x)^{1/2}}\log x\right).\end{aligned} \tag{101}$$

To complete the proof observe that $O(xe^{-c(\log x)^{1/2}}\log x)=O(xe^{-c_0(\log x)^{1/2}})$ for some constants $c > 0$, and $c_0 > 0$. ■

This proof shows that the estimate $M(x)=\sum_{n\le x}\mu(n)=O(xe^{-c(\log x)^{1/2}})$ of the summatory Mobius function implies the prime number theorem. Other proofs are given in [MV, p. 179], [NW, p. 232], [TM, p. 167], and the literature. The best error term known, obtained by complex analysis, is

$$\pi(x)=li(x)+O\left(xe^{-c(\log x)^{3/5}(\log\log x)^{-1/5}}\right), \tag{102}$$

see [FD], and the best error term known, obtained by elementary methods, is

$$\pi(x)=li(x)+O\left(xe^{-c(\log x)^{1/10}}\right), \tag{103}$$

see [DD].

The basic principle on the oscillatory nature of the function ψ(*x*) − *x* is well known in the literature, [EL, p. 186], [MV, p. 464], [NW, p. 226]. A refined version is stated below.

The big omega symbol is defined as follows:

$$\text{(i)}\quad f(x)=\Omega_+(g(x)) \iff \limsup_{x\to\infty}\frac{f(x)}{g(x)}>0, \qquad \text{(ii)}\quad f(x)=\Omega_-(g(x)) \iff \limsup_{x\to\infty}\frac{f(x)}{g(x)}<0, \tag{104}$$

see [EL, p. 186], MV, p. 4], [WK] et cetera.

***Theorem* 7.3.** Let $x \ge 1$ be a large real number. Then

$$\text{(i)}\quad \pi(x)=li(x)+\Omega_\pm\left(x^{1/2}\log\log\log x/\log x\right), \tag{105}$$

$$\text{(ii)}\quad \pi(x)=li(x)+\Omega_\pm\left(x^{1/2}\log\log\log x/\log x\right),$$

where $c > 0$ is an absolute constant.

The proofs are given in [EL, p. 191], [IV, 307], [MV, p. 478], [NW, p. 322], and similar reference.

***Theorem* 7.4.** Let $x \ge 1$ be a large real number. Let $a \ge 1$, $q \ge 1$ be a pair of integers such that gcd($a$, $q$) = 1, and let $q = O((\log x)^B)$, $B > 0$ constant. Then

$$\text{(i)}\quad \pi(x,a,q)=\frac{li(x)}{\varphi(q)}+O\left(xe^{-c(\log x)^{1/2}}\right), \tag{106}$$

$$\text{(ii)}\quad \psi(x,a,q)=\frac{x}{\varphi(q)}+O\left(xe^{-c(\log x)^{1/2}}\right),$$

where $c > 0$ is an absolute constant.

Confer [MV, p. 382], [NW], [TM, p. 255], and similar literature.

The range of $q \ge 1$ can be extended to $q \le x$, but the asymptotic formula has a weaker error term, namely,

$$\pi(x,a,q) = \frac{x}{\varphi(q)\log x} + O(\frac{x}{\log^B x}), \tag{107}$$

where $B > 0$, see [EL].

## 8. Elementary Results For Quadratic Arithmetic Progressions

The basic results required to settle the existence of infinitely many primes in certain quadratic arithmetic progressions are assembled here.

### 8.1 Powers Sums Over Quadratic Arithmetic Progressions

Let $a$, and $q \geq 1$ be integers, $\gcd(a, q) = 1$. The element $a \in \mathbb{Z}_q = \{ 0, 1, 2, \ldots, q - 1 \}$ is called a *quadratic residue* the congruence $n^2 \equiv a \bmod q$ has a solution. Otherwise $a$ is a *quadratic nonresidue*. The quadratic symbol is defined by

$$\left(\frac{x}{p}\right) = \begin{cases} 1 & \text{if } x \text{ is a quadratic residue,} \\ -1 & \text{if } x \text{ is a quadratic nonresidue,} \\ 0 & \text{if } x \ \gcd(x,q) \neq 1. \end{cases} \tag{108}$$

Often, the quadratic aspect of an element in the residue number system can be determined via the quadratic reciprocity law and related formulas:

(i) $\left(\frac{-1}{p}\right) = (-1)^{(p-1)/2}$, (ii) $\left(\frac{2}{p}\right) = (-1)^{(p^2-1)/8}$,

(iii) $\left(\frac{p}{q}\right)\left(\frac{q}{p}\right) = (-1)^{(p-1)(q-1)/4}$. (109)

For a modulo $q \geq 1$, and a quadratic residue $a \geq 1$, with $\gcd(a, q) = 1$, the congruence $n^2 \equiv a \bmod q$ has $2^W \geq 2$ solutions, where $W = \omega(q) + r$, and $r = 0, 1, 2$ according as 4 does not divide $q$, $4 \parallel q$ or $8 \mid q$, see [LV, p. 65]. The function $\omega(n) = \#\{ p : p \mid n \}$ tallies the number of prime divisors in $q$.

These concepts come into play in the analysis of power sums over quadratic arithmetic progressions.

***Lemma* 8.1.** Let $a \geq 1$, and $q \geq 1$ be integers, $\gcd(a, q) = 1$. Assume that $a$ is a quadratic residue modulo $q$. Then

(i) $\displaystyle\sum_{n \leq x,\ n^2 \equiv a \bmod q} n^2 = \frac{2^W}{3q^2} x^3 + O\left(x^2\right)$, (ii) $\displaystyle\sum_{n \leq x,\ n^2 \equiv a \bmod q} n^2 \geq \frac{2^W}{3q^2} x^3 + O\left(\frac{1}{q} x^2\right)$, (110)

where $2^W \geq 2^{\omega(q)}$ with $\omega(n) = \#\{ p : p \mid n \}$, for all sufficiently large real numbers $x \geq 1$.

Proof (ii): For each fixed $q \leq x$ the congruence $n^2 \equiv a \bmod q$ has at least $2^W \geq 2^{\omega(q)} \geq 2$ solutions if and only if $a \geq 1$ is a quadratic residue $b^2 \equiv a \bmod q$ with $0 < b < q$. For each such solution $b < q$ the sequence of integers in the arithmetic progression $n = qm + b$, with $0 \leq m \leq (x - b)/q$, are also solutions of the congruence. Ergo, expanding, substituting Lemma 8.3, and simplifying return

$$\sum_{n \le x,\ n^2 \equiv a \bmod q} n^2 = \sum_{n \le x,\ n \equiv b \bmod q} n^2$$

$$= q^2 \sum_{m \le (x-b)/q} m^2 + 2bq \sum_{m \le (x-b)/q} m + b^2 \sum_{m \le (x-b)/q} 1$$

$$\ge q^2 \left( \frac{1}{3q} \left( \frac{x-b}{q} \right)^3 + O\left( \frac{1}{q} (\frac{x-b}{q})^2 \right) \right) + 2bq \left( \frac{1}{2q} \left( \frac{x-b}{q} \right)^2 + O\left( \frac{x-b}{q} \right) \right) + b^2 \left( \frac{x-b}{q} + O(1) \right) \tag{111}$$

$$\ge \frac{2^W}{3q^2} x^3 + O\left( \frac{1}{q} x^2 \right),$$

where the value $x^3/3q^2 + O(x^2/q) \ge 0$ for all triples $0 < a < q \le x$. Lastly, the $2^W$ distinct roots $b \in \{ x : x^2 \equiv a \bmod q \ne b \bmod q \}$ generate $2^W$ distinct arithmetic progressions, so the value $x^3/3q^2 + O(x^2/q)$ is $2^W$ fold. ■

The case of interest in Theorem 4.1 calls for $a = -1$, and $q = d$, so the prime divisors $p$ of $q$ satisfy the congruence $p \equiv 1 \bmod 4$, and $2^W \ge 2^{\omega(d)} \ge 2$.

## 8.2 Twisted Finite Sums Over Quadratic Arithmetic Progressions

***Lemma* 8.2.** Let $s \ge 1$ be a real number, let $\mu(n)$ be the Mobius function, and let $x \ge 1$ be a sufficiently large real number. Then

$$-\sum_{n \le x,\ w^2 \equiv -1 \bmod n} \frac{\mu(n)\log n}{n^s} \ge c_1 + O\left( x^{1-s} e^{-c(\log x)^{1/2}} \log x \right), \tag{112}$$

where $c_1 > 0$, and $c > 0$ are constants.

Proof: First compare the two twisted finite sums:

$$-\sum_{\substack{n \le x, \\ w^2 \equiv -1 \bmod n}} \frac{\mu(n)\log n}{n} = -\sum_{\substack{n \le x, \\ p \mid n \Rightarrow p \equiv 1 \bmod 4}} \frac{\mu(n)\log n}{n}$$

$$= \sum_{\substack{p \le x, \\ p \equiv 1 \bmod 4}} \frac{\log p}{p} - \sum_{\substack{p, q \le x, \\ p, q \equiv 1 \bmod 4}} \frac{\log pq}{pq} + \sum_{\substack{p, q, r \le x, \\ p, q, r \equiv 1 \bmod 4}} \frac{\log pqr}{pqr} - \cdots \tag{113}$$

and

$$-\sum_{\substack{n \le x, \\ n \equiv 1 \bmod 4}} \frac{\mu(n)\log n}{n} = \sum_{\substack{p \le x, \\ p \equiv 1 \bmod 4}} \frac{\log p}{p} - \sum_{\substack{pq \le x, \\ p, q \equiv 1 \bmod 4}} \frac{\log pq}{pq} - \sum_{\substack{pq \le x, \\ p, q \equiv 3 \bmod 4}} \frac{\log pq}{pq} + \sum_{\substack{pqr \le x, \\ p, q, r \equiv 1 \bmod 4}} \frac{\log pqr}{pqr} + \cdots . \tag{114}$$

Clearly, the difference

$$\left( - \sum_{\substack{n \le x, \\ w^2 \equiv -1 \bmod n}} \frac{\mu(n)\log n}{n^s} \right) - \left( - \sum_{\substack{n \le x, \\ n \equiv 1 \bmod 4}} \frac{\mu(n)\log n}{n^s} \right) = \sum_{\substack{pq \le x, \\ p,q \equiv 3 \bmod 4}} \frac{\log pq}{p^s q^s} - \sum_{\substack{pqr \le x, \\ p,q \equiv 3 \bmod 4,\ r \equiv 1 \bmod 4}} \frac{\log pqr}{p^s q^s r^s} + \cdots > 0 \tag{115}$$

is a nonnegative real number. Therefore

$$\begin{aligned} - \sum_{n \le x,\ w^2 \equiv -1 \bmod n} \frac{\mu(n)\log n}{n^s} &= - \sum_{n \le x,\ p \mid n \Rightarrow p \equiv 1 \bmod 4} \frac{\mu(n)\log n}{n^s} \\ &\ge - \sum_{n \le x,\ n \equiv 1 \bmod 4} \frac{\mu(n)\log n}{n^s} \\ &= \kappa_1 + O\left( x^{1-s} e^{-c(\log x)^{1/2}} \log x \right), \end{aligned} \tag{116}$$

where $\kappa_1 = -L'(s,\chi_0)L(s,\chi_0)^{-2}$ is a constant, this follows from Lemma 8.4. ■

In terms of L-function, the constant has the general form

$$- \sum_{n \equiv 1 \bmod 4} \frac{\mu(n)\log n}{n^s} = \frac{d}{ds} \frac{1}{L(s,\chi_0)} = - \frac{L'(s,\chi_0)}{L(s,\chi_0)^2}, \tag{117}$$

where $L(s,\chi_0) = \sum_{\ge 1} \chi_0(n) n^{-s}$ is a Dirichlet series, and $\chi_0 = 1$ if gcd(n, q) = 1, else $\chi_0 = 0$, is the principal character modulo $q$, see [DL, 25.15.1] for some details on L-function.

### 8.3 Powers Sums Over Arithmetic Progressions

Let $\mathbb{N} = \{ 0, 1, 2, 3, \ldots \}$ be the set of nonnegative integers, and let $q\mathbb{N} + a = \{ qn + a : n \ge 0 \}$ be the arithmetic progression defined be a pair of integers $a \ge 0$, and $q \ge 1$.

The standard formulas for the sums of powers over the integers are

$$\sum_{n \le x} n = \frac{1}{2}\left(x^2 + x\right), \quad \sum_{n \le x} n^2 = \frac{1}{6}\left(2x^3 + 3x^2 + x\right), \quad \sum_{n \le x} n^2 = \frac{1}{6}\left(2x^3 + 3x^2 + x\right),$$

$$\sum_{n \le x} n^3 = \frac{1}{4}\left(x^4 + 2x^3 + x^2\right), \ldots . \tag{118}$$

These are well know, and there are numerous formulas, and technique for evaluating these finite sums. The sums of powers over arithmetic progressions is one of the possible generalizations of these finite sums.

***Lemma* 8.3.** Let $a \ge 0$, and $q \ge 1$ be fixed integers. Let $x \ge 1$ be a sufficiently large real number. Then

$$\text{(i)} \sum_{n \le x,\ n \equiv a \bmod q} n = \frac{1}{2q} x^2 + O(x). \qquad \text{(ii)} \sum_{n \le x,\ n \equiv a \bmod q} n \ge \frac{1}{2q} x^2 + O\left(\frac{1}{q} x\right). \tag{119}$$

$$\text{(iii)} \sum_{n \le x,\ n \equiv a \bmod q} n^2 = \frac{1}{3q^2} x^3 + O\left(x^2\right). \qquad \text{(iv)} \sum_{n \le x,\ n \equiv a \bmod q} n^2 \ge \frac{1}{3q^2} x^3 + O\left(\frac{1}{q} x^2\right).$$

Proof (i): The integers in the linear arithmetic progression are of the form $n = qm + a$, with $0 \le m \le (x - a)/q$. Inserting this into the finite sum produces

$$\begin{aligned}\sum_{n \le x,\, n \equiv a \bmod q} n &= q \sum_{m \le (x-a)/q} m \;+\; a \sum_{m \le (x-a)/q} 1 \\ &= \frac{q}{2}\left(\frac{x-a}{q}\right)\left(\frac{x-a}{q}+1\right) \;+\; a\left(\frac{x-a}{q}\right).\end{aligned} \tag{120}$$

This last expression simplifies into the expression given above. Likewise for (ii), expand it and substitute the standard sums of powers:

$$\begin{aligned}\sum_{n \le x,\; n \equiv a \bmod q} n^2 &= q^2 \sum_{m \le (x-a)/q} m^2 \;+\; 2aq \sum_{m \le (x-a)/q} m \;+\; a^2 \sum_{m \le (x-a)/q} 1 \\ &= \frac{q^2}{6}\left(2\left(\frac{x-a}{q}\right)^3 + 3\left(\frac{x-a}{q}\right)^2 + \frac{x-a}{q}\right) \;+\; 2bq\left(\frac{1}{2}\left(\frac{x-a}{q}\right)\left(\frac{x-a}{q}+1\right)\right) \;+\; a^2\left(\frac{x-a}{q}\right) \\ &= \frac{1}{3q^2}x^3 + O\left(x^2\right),\end{aligned} \tag{121}$$

where the standard formulas are given in (118). The last expression simplifies into the expression (iii) given above. To obtain the lower estimate in (iv), examine the explicit expansion

$$\sum_{n \le x,\; n \equiv a \bmod q} n^2 = \sum_{n \le x,\; n \equiv b \bmod q} n^2 = \frac{1}{3q}x^3 - \frac{a}{3q}x^2 + \cdots \ge \frac{1}{3q}x^3 + O(\frac{1}{q}x^2), \tag{122}$$

and take the smallest possible value for any $q \le x$. ■

# 9. Application of Germain Primes

There are several applications of Germain primes to other problems in the Mathematical sciences. These applications ranges from the theory of Diophantine equations, to primality testing. An application will be sketched in this Section.

## 9.1. Application To Primitive Roots

Let $G$ be a cyclic group of order $N = \#G$. The important primitive root tests, procedures (i) and (ii) below, of primitive root in $G$ are directly dependent on the prime decomposition of the integer $\varphi(N) = N\prod_{p|N}(1-1/p)$. This dependence on the prime decomposition of the integer $N \ge 1$ induces a complicated mechanism on the distribution of primitive roots modulo $N$.

(i) The Lehmer primitive root test: Let $r > 1$ be a primitive root modulo $N$. Then

$$r^{\varphi(N)/p} \not\equiv 1 \bmod N \tag{123}$$

for every prime divisor $p$ of $N$, confer [CE, p. 25], and [CP].

(ii) The Characteristic function primitive root test: An element $r \in \mathbb{Z}_N = \{$ residues modulo $N \}$ is a primitive root modulo $N$ if and only if

$$\frac{\varphi(N)}{N}\sum_{d \mid N}\frac{\mu(d)}{\varphi(d)}\sum_{ord(\chi)=d}\chi(r) > 0, \tag{124}$$

where $\text{ord}(\chi) = \min \{ d \geq 1 : \chi_d = 1 \}$ denotes the order of the character $\chi$ modulo $N$, confer [MO p. 17].

These observations shows that the complexity of computing a primitive root in a cyclic group depends on the parameters

$$\omega(N) \geq 1, \qquad \varphi(N)-(N-1)/2, \tag{125}$$

and other intrinsic parameters of the integer $N \geq 1$.

For arbitrary integers $N$, these primitive root tests are difficult to verify. However, there are sequences of special prime numbers, which have simple structures and simple primitive root verification procedures. For example, the sequence of Fermat primes $F_n = 2^{2^n} + 1$, $n \geq 0$, has the simplest primitive root verification procedure, called the Pepin test. In fact, any odd number is a primitive root modulo $F_n$ since $\varphi(F_n - 1) - (F_n - 1)/2 = 0$ is minimal. Other sequences of prime numbers, which has a quite simple primitive root verification procedures, are the sequence of prime numbers $p = 2q + 1$, $p = 4q + 1$, et cetera.

***Theorem* 9.1.** There are infinitely many primitive roots $g$ modulo $p$ in following cases.
(i) $g = 2$, and $p = 2q + 1 = 8n + 3$ is a Germain prime.
(ii) $g = -2$, and $p = 2q + 1 = 8n + 7$ is a Germain prime.
(iii) $g = 2$, and $p = 4q + 1 = 8n + 5$ is a Germain prime.

Proof (i): Since integer $p - 1 = 2q$ has just two prime divisors 2 and $q$, the primitive root verification procedure states that 2 is a primitive root modulo $p = 2q + 1 = 8n + 3$ if and only if

$$\text{(i)} \quad 2^{(p-1)/2} \not\equiv 1 \bmod p, \qquad \text{and} \qquad \text{(ii)} \quad 2^{(p-1)/q} \not\equiv 1 \bmod p, \tag{126}$$

see [CE, p.25], [CP], [LV, 68]. Using the quadratic reciprocity law, it easily follows that both conditions

$$2^{(p-1)/2} \equiv \left(\frac{2}{p}\right) \equiv (-1)^{(p^2-1)/8} \equiv -1 \bmod p, \qquad \text{and} \qquad 2^{(p-1)/q} \equiv 2^2 \bmod p \tag{127}$$

are satisfied modulo $p = 2q + 1 = 8m + 3$. Ergo, the integer 2 is a primitive root modulo $p = 2q + 1 = 8n + 3$. Moreover, by Theorem 3.1, it quickly follows that the lower bound of the Germain prime counting function on arithmetic progression $\{ p = 8n + 3 : n \geq 1 \}$ satisfies

$$\pi_G(x,a,q) \geq \frac{c_0}{\varphi(q)}\frac{x}{\log^2 x} + o\left(\frac{x}{\log^2 x}\right), \tag{128}$$

where $c_0 > 0$ is a constant, shows that there are infinitely many Germain primes in the arithmetic progression $p = 2q + 1 = 8m + 3$. The proofs of (ii) and (iii) are identical, but use other residue classes $p = 2q + 1 = 8m + 7$, and $p = 2q + 1 = 8m + 5$. ■

**9.2. The Primitive Roots Problem**

The only known case of a fixed primitive $g$ mod $N$ of an infinite sequence of integers $N$ is $g \equiv (-1)^r 5^s \bmod N = 2^m$, $m \geq 1$, and ($r \geq 0$, $s \geq 0$). But, the pair $< -1, 5 >$ generates the multiplicative group of $\mathbb{Z}_N$ using two elements, not a single element the multiplicative group of $\mathbb{Z}_N$.

The primitive root problem, best known as Artin conjecture, states the followings.

***Primitive Roots Problem*** (Artin conjecture) Let $a \in \mathbb{Z}$ be an integer, $|a| \geq 1$, and $a \neq b^2 \not\equiv 1 \bmod 4$. Then, there are

$$\pi_a(x) = C_a \frac{x}{\log x} + o\left(\frac{x \log\log x}{\log^2 x}\right) \tag{129}$$

prime numbers $p \leq x$ such that $a$ is a primitive root modulo $p$. The constant is defined by

$$C_a = \prod_{p \geq 2}\left(1 - \frac{1}{p(p-1)}\right) = .373956..., \tag{130}$$

see [MU], [MO], et alii, for surveys and background information.

For integers $a \in \mathbb{Z}$ of the form $a = b^2c \equiv 1 \bmod 4$, the constant $C_a$ is readjusted to reflect the other primes interdependence, see [BS, p.254], [LN, p. 2] and similar source.

The density of a sequence of Germain primes $p = 2q + 1$ is smaller than the conjecture density of primitive roots by a factor of $\log x$. Thus, Theorem 9.1 is an ad hoc method of constructing examples of the Artin conjecture. Even the union of log $0x$ distinct sequences of Germain or other primes does not seem to be sufficient to settle the primitive root problem since the asymptotic expressions do not seem to agree. The complete solution of the primitive root conjecture is a much more difficult problem.